\documentclass{article}
\usepackage{amsmath,amsthm, amsfonts, amssymb, amscd, enumerate,mathrsfs,authblk}
\usepackage{booktabs} 
\usepackage{array} 
\usepackage{paralist} 
\usepackage{subfig} 
\usepackage{fancyhdr} 
\usepackage {epsfig}
\usepackage{graphicx}
\usepackage{subcaption}
\usepackage {url}
\usepackage {textcomp}
\usepackage{verbatim}
\usepackage[all]{xy}
\usepackage{bm}
\usepackage[export]{adjustbox}
\usepackage{wrapfig}
\usepackage{enumitem}
\usepackage{placeins}
\newcommand{\rr}{\mathbb{R}}

\newcommand{\mcL}{\mathcal{L}}
\newcommand{\mcX}{\mathcal{X}}

\newcommand{\mcI}{\mathcal{I}}
\newcommand{\mcD}{\mathcal{D}}
\newcommand{\mcB}{\mathcal{B}}

\newcommand{\norm}[1]{\left\Vert{#1}\right\Vert}

\theoremstyle{plain}
\theoremstyle{definition}

\title{The RBF Collocation Method to Design a Digital Twin for Coffee Percolation}
\author[1]{Nadaniela Egidi\footnote{nadaniela.egidi@unicam.it}}
\author[2]{Lauro Fioretti\footnote{lauro.fioretti@simonelligroup.it}}
\author[1]{Josephin Giacomini\footnote{josephin.giacomini@unicam.it}}
\author[1]{Pierluigi Maponi\footnote{pierluigi.maponi@unicam.it}}
\author[1]{Gianluca Pacini\footnote{gianluca.pacini@unicam.it}}

\affil[1]{School of Science and Technology - Mathematics Division, University of Camerino, Via Madonna delle Carceri 9, 62032 Camerino,(MC), Italy}

\affil[2]{Simonelli Group}

\begin{document}

\maketitle

\begin{abstract}
    Espresso coffee extraction is a complex physico-chemical process and can be modeled through a system of coupled partial differential equations. We present a numerical solution based on a meshless Collocation Method using Radial Basis Functions and Kansa’s approach, which reveals to be accurate and robust in comparison to a reference numerical solution provided by a well-known simulation software.
\end{abstract}

\section{Introduction}
Coffee is one of the most widely consumed beverages in the world after water and tea. It is estimated that over 2 billion cups are consumed every day. For this reason, many models that simulate the extraction process of espresso coffee are studied. These models consist of a physico-chemical process: hot water at a given pressure pours into the basket and hits the tamped coffee powder within the basket, it flows into the powder, passing through the void spaces between the coffee grains of the coffee pod. The water also removes a certain amount of fine particles from the ground coffee and transports them downward. This extraction process can be described as a fluid-dynamics process that is a well-known percolation process. 

From the fluid-dynamic point of view, these components are involved: the dynamics of the fluid, that is described by the Darcy law~\cite{Darcy}; the behavior of the pressure, that is ruled by the Richards' equation~\cite{MGFaCT},~\cite{bear2013dynamics},~\cite{pinder2008essentials},~\cite{verruijt2017introduction}; the transport, and the dissolution of the chemical species~\cite{CiPM},~\cite{HoPM},~\cite{FEFLOW}; the transfer of heat between the solid medium and the fluid, described by the heat equation~\cite{CiPM}.

There exist various numerical methods for the approximate solution of Partial Differential Equations (PDEs) with initial and boundary conditions~\cite{NAoPDE}.
In this article, a model describing the coffee percolation process has been presented and solved numerically. In particular, the Collocation Method joined with a Radial Basis Functions (RBFs) approximation is used to numerically solve the considered problem with Kansa’s approach~\cite{MAMM},~\cite{SDA}.

In Section \ref{sec_CPM}, the percolation model is described, and the initial and boundary conditions are defined.
In Section \ref{sec_KM}, the numerical methods we used to solve the model are explained.
In Section \ref{sec_PA}, describes how the model and domain were discretized.
In Section \ref{sec_NR}, the results obtained are shown.

\section{Coffee Percolation Model}
\label{sec_CPM}
\begin{figure}
    \centering
    \includegraphics[width=0.7\textwidth]{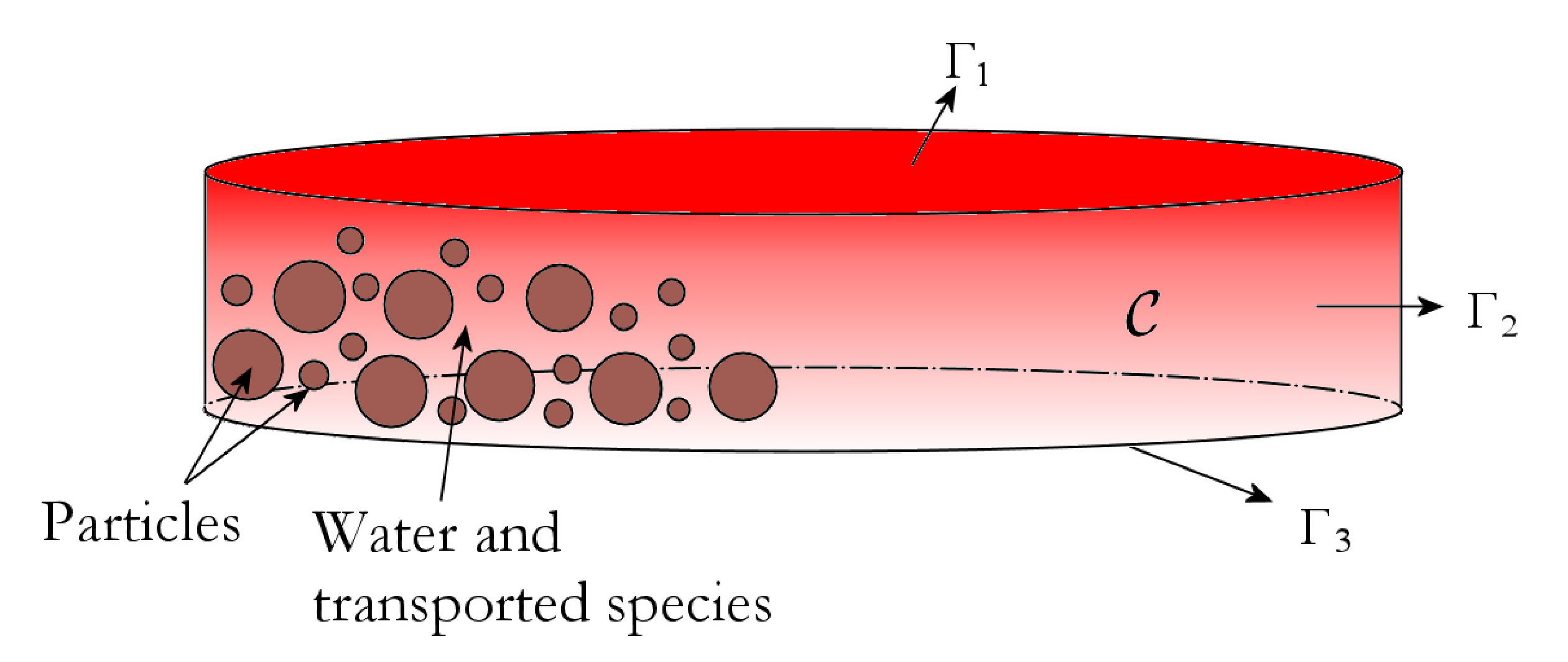}
    \caption{Domain $\mathcal{D}= \mathcal{C} \cup\Gamma_1 \cup\Gamma_2 \cup\Gamma_3$ of the Percolation Model, i.e., schematization of the coffee pod.}
    \label{CofPod}
\end{figure}

The coffee percolation model describes the behavior, in the coffee pod represented in Fugure \ref{CofPod}, of the following quantities: the hydraulic head $h(t,\bm x)$, at time $t\ge 0$ and point $\bm x=(x_1,x_2,x_3)\in\mathcal{D}\subset\rr^3,$ which is a quantity closely related to the pressure $p(t,\bm x)$ by the relation $h=p/\rho_0 g +x_3$, where $\rho_0$ is the fluid mass density of reference and $g$ is the gravitational acceleration; the temperature $T(t,\bm x)$ of the coffee pod; the transport of $N_s$ chemical species, in particular with $C_k(t,\bm x)$, $k=1,2,\dots,N_s,$ we denote the concentration of the $k$th species in the liquid, and we will call it the liquid concentration, and their dissolution, we denote with $C_k^s(t,\bm x)$ the solid concentration of the chemical species $k$.
%The  model describe the behavior of the hydraulic head $h(t,\bm x)$, at time $t$ and point $\bm x=(x_1,x_2,x_3)\in\rr^3,$ which is a quantity closely related to the pressure $p(t,\bm x)$ by the relation $h=p/\rho_0 g +x_3$, where $\rho_0$ is the fluid mass density of reference and $g$ is the gravitational acceleration; the behavior of the temperature $T(t,\bm x)$ of the coffee pod; the transport of chemical species $C_k$, which we will call liquid, and the dissolution of the chemical species $C_k^s$, which we will call solid.

To build the model, we made the following assumptions:
\begin{enumerate}
        \item the porous medium has the same physical properties in all directions, i.e., it is isotropic;
        \item the coffee powder is ground with different sizes of coffee grains, but we assume that their physical and chemical properties are the same in each grain, independently of their size and location, i.e., the porous medium is homogeneous;
        \item the percolation problem lasts about 25 seconds, but the coffee espresso comes out after about 5 seconds, which is necessary for the imbibition of the coffee pod. After the imbibition, there is no gaseous phase in the porous medium, i.e., we consider the porous medium saturated, and a local thermal balance between the coffee powder and water;
        \item the flow is smooth and ordered which fluid particles move in parallel layers with minimal mixing between them, i.e., the flux is laminar;
        \item we neglect the transport of the fine particles, and therefore the creation of the compact layer, i.e., we consider the porous medium fixed;
        \item there are no internal sinks or sources of water, heat, or chemical substances in the coffee pod.
\end{enumerate}

Hence, the equations that describe the coffee percolation are:
\begin{equation}
\begin{cases}
    S_0\frac{\partial h}{\partial t}+\nabla\cdot\bm{q} = 0,\\
    \bm{q}=-\bm{K}f_\mu\cdot\left( \nabla h+\chi\bm{e}\right),\\
    \epsilon\frac{\partial C_k}{\partial t}+\bm{q}\cdot\nabla C_k+\nabla\cdot\bm{j_k}=R_k, \qquad k=1,\dots,N_{s},\\
    \epsilon_s\frac{\partial C_k^s}{\partial t}=R_k^s, \qquad k=1,\dots,N_{s},\\
    \left(\epsilon\rho c+\epsilon_s\rho^s c^s \right)\frac{\partial T}{\partial t}+\rho c \bm{q}\cdot \nabla T-\nabla \cdot \left( \bm{\Lambda} \cdot \nabla T \right)=0,
\end{cases}\label{model}
\end{equation}
where $S_0$ is the specific storage coefficient, the column vector $\bm{q}(t,\bm x)\in\rr^3$ is the Darcy flux, $\bm{K}\in\rr^{3\times 3}$ is the hydraulic conductivity tensor, $f_\mu$ is the viscosity relation function, $\chi$ is the buoyancy coefficient, $\bm{e}=(0,0,1)^T\in \rr^3$ (the apex $T$ denotes the transposte), $\epsilon$ is the porosity of the medium, $\bm{j_k}=\bm{D_k}\cdot\nabla C_k$ is the diffusive flux tensor where $\bm{D_k}\in\rr^{3\times 3}$ is the hydrodynamic dispersion tensor, consisting of the sum of molecular diffusion and mechanical diffusion, and is given by
\[
\bm{D_k}=\left(\epsilon D_k+\beta_T^k\norm{\bm{q}}\right)I+ \left(\beta_L^k + \beta_T^k\right)\frac{\bm{q}\otimes \bm{q}}{\norm{\bm{q}}},
\]
where $\beta_L^k,\beta_T^k$ are the transverse and the longitudinal dispersion coefficients, respectively, $D_k$ is the molecular diffusion coefficient, $I\in\rr^{3\times 3}$ is the identity matrix, and the symbol $\otimes$ denotes the tensor product.

Moreover, $\epsilon_s$ is the solid volume fraction, and it satisfies $\epsilon_s=1-\epsilon$, $\rho^s,\rho$ are the density of solid and fluid, $c^s,c$ are the specific heat of the solid and the fluid (so $\rho c$ is the fluid volumetric heat capacity and $\rho^s c^s$ is the solid volumetric heat capacity), the total reaction rate terms are defined as:
\begin{equation}
R_k=\alpha_k\epsilon_sC_k^s, \quad R_k^s=-\alpha_k\epsilon_s C_k^s, \quad \epsilon_s=1-\epsilon,
\label{RkRks}
\end{equation}

where the coefficients $\alpha_k$ are function of incoming water pressure $p_{z0}$ and incoming water temperature $T_{z0}$:
\[
\alpha_k =A_0 + aT_{z0} + bp_{z0} + cT_{z0}^2 + dp_{z0}^2 +fT_{z0}p_{z0} + lT_{z0}^2p_{z0} + mT_{z0}p_{z0}^2
\]
where the coefficients $A_0,a,b,c,d,f,l,m$ depend on the chemical species, the granulometry of coffee powder, and the coffee blends. In the end, $\bm{\Lambda}\in \rr^{3 \times 3}$ is the dispersion tensor and for an isotropic medium has a similar shape to that of the hydrodynamic dispersion tensor $\bm{D_k}$.

By considering Figure \ref{CofPod}, let be ${\bm n}$ the normal vector exiting a closed surface $\partial \mathcal{D}$, $r=\norm{\bm x}$ the Eucledean norm,  $h_{z0},\Phi_h,h_C,C_{kC},C_{k0}^s,T_{z0},T_{0}\in\rr$ given constants, and $p_0(\cdot)$ a given function on $\rr.$
The boundary and initial conditions are the following:
 \begin{equation}
\begin{cases}
    h=h_{z0}, & \text{on } \Gamma_1,\ t>0,\\
    \frac{\partial h}{\partial r}=0, & \text{on } \Gamma_2,\ t>0,\\
    \bm{q}\cdot \bm{n}=-\Phi_h\text{min}\{h_C-h,0\}, & \text{on } \Gamma_3,\ t>0,\\
    p=p_0(x_3), & \text{on } \mathcal{D},\ t=0, \label{RicDar}
\end{cases}
\end{equation}
\begin{equation}
\begin{cases}
    \nabla C_k\cdot \bm{n}=0, \qquad & \text{on } \Gamma_1,\Gamma_2,\ t>0,\\
    -\left(\bm{D_k}\cdot\nabla C_k\right)\cdot \bm{n}=-\Phi_k\text{min}\{C_{kC}-C_k,0\}, & \text{on } \Gamma_3, t>0,\\
    C_k=0, & \text{on } \mathcal{D},\ t=0,
\end{cases}
\end{equation}
\begin{equation}
    C_k^s=C_{k0}^s, \qquad \text{on } \mathcal{D},\ t>0,\label{ICDiss}
\end{equation}
\begin{equation}
\begin{cases}
    T=T_{z0}, & \text{on } \Gamma_1,\ t>0,\\
    \nabla T\cdot \bm{n}=0, & \text{on } \Gamma_2,\Gamma_3,\ t>0,\\
    T=T_0, & \text{on } \mathcal{D},\ t=0.
    \end{cases}
\end{equation}

For more details about how the variables of the model are chosen, see~\cite{Art_perc}.

We note that the mass dissolution equations, the fourths in \eqref{model}, with initial condition \eqref{ICDiss}, and relations \eqref{RkRks} have the following analytical solutions
\begin{equation}
    C_k^s(t,\bm x)=C_{k0}^se^{-t\alpha_k},\qquad t\ge 0,\ \bm x\in \mcD.\label{AnSolDiss}
\end{equation}

In the following, for the heat equation, the fifth in \eqref{model}, we consider the following approximation: 
\[
\bm{\Lambda} =\begin{bmatrix}
    \lambda & 0 & 0\\
    0 & \lambda & 0\\
    0 & 0 & \lambda
\end{bmatrix}
=\lambda\bm{I}
\Longrightarrow 
\nabla \cdot \left( \bm{\Lambda} \cdot \nabla T \right) = \lambda\Delta T,
\]
where $\lambda\in\rr$. Moreover, in the mass transport equations, the fourths in \eqref{model}, we consider the Darcy flux $\bm{q}$ constant and equal to $(0,0,q_0)$, where $q_0$ is found by solving the Richards' equation, the first two in \eqref{model}.

\section{Kansa's Method}
\label{sec_KM}
A function $\Phi:\rr^s\rightarrow \rr$ is called radial if there exists a univariate function $\phi:\left[0,\infty\right]\rightarrow \rr$ such that
\[
\Phi(\bm{x})=\phi(\norm{\bm{x}}), \qquad {\bm{x}\in\rr^s},
\]
and $\norm{\cdot}$ is the Euclidean norm on $\rr^s$.
This definition says that for a radial function $\Phi$
\[
\norm{\bm x_1} = \norm{\bm x_2} \Longrightarrow \Phi(\bm x_1
)=\Phi(\bm x_2), \qquad \bm x_1,\bm x_2\in\rr^s;
\]
thus, $\Phi$ is radially symmetric about its center.

We use the following notations: $\Omega$ denotes an open subset of $\rr^s$, $\overline\Omega$ is its closure, $\partial \Omega$ is its boundary, $\partial^* \Omega\subset \partial \Omega $,  $\rr^+_0$ is the set of non negative real numbers,  $\rr^+$ is the set of positive real numbers.

Given $\mathcal{X}=\{\bm x_1,\bm x_2,\dots,\bm x_N\}\subset\overline\Omega,$  a set of points we choose an approximation of a function $u(t,\bm{x})$, $\bm{x}\in\Omega\subseteq\rr^s$ and $t\in\rr^+_0$, in the family $\mathcal{F}\left(\mathcal{X},M,\phi\right)$ of functions having the following form  
\begin{equation}
\hat u(t,\bm{x})=\displaystyle\sum_{i=1}^Nu_i(t)\phi(\norm{\bm{x}-\bm x_i})+\displaystyle\sum_{i=N+1}^{N+M}u_i(t)P_{i-N}(\bm{x}),\qquad (t,\bm x)\in \rr_0^+\times \overline \Omega\label{hatu}
\end{equation}
where $M>0,$ $u_i(t)$ are differentiable functions for $i=1,2,\dots, N+M,$ $\phi$ is an RBF, and $P_k$ is a polynomial of degree at most $M$, for a comprehensive discussion see~\cite{MAMM}.

We consider PDEs of the form 
\begin{equation}
     \frac{\partial u(t,\bm{x})}{\partial t} +\mathcal{L}u(t,\bm{x}) = f(t,\bm{x}),\qquad (t,\bm x)\in \rr^+\times, \Omega\label{PDE}
\end{equation}
with boundary condition
\begin{equation}
     \mathcal{B}u(t,\bm{x}) = g(t,\bm{x}), \qquad (t,\bm x)\in \rr^+\times \partial^* \Omega\label{BPDE}
\end{equation}
where $\mcL$ and $\mathcal{B}$ are differential operators with respect to $\bm{x}$.

Hence, by using Kansa's collocation method with the set of collocation nodes equal to $\mcX,$ we want to compute column vectorial function  $\bm u:\rr^+_0\to\rr^{N+M}$ having components $u_i(t)$, $i=1,2,\dots, N+M,$ $t\in\rr_0^+,$ such that the approximation $\hat u$ of $u$ satisfies \eqref{PDE} on $\rr^+\times (\Omega\cap \mcX) $ and   \eqref{BPDE} on $\rr^+\times (\partial^*\Omega\cap \mcX). $
Let be 
$$\mcI=\{1,2,\dots,N\},\quad \mcI_0=\{i\in\ \mcI\ :\ \bm x_i\in \Omega\}, \quad \mcI_\partial=\{i\in\ \mcI\ :\ \bm x_i\in \partial ^*\Omega\},$$
we suppose that $\mcI=\mcI_0\cup \mcI_\partial$ and $\mcI_0\cap \mcI_\partial=\emptyset.$ Then we obtain the following system of Differential Algebraic Equations (DAE):

for $j\in \mcI_0$
\begin{equation}
     \displaystyle\sum_{i=1}^N\phi(\norm{\bm{x}_j-\bm{x}_i})\dot{u}_i(t)+\displaystyle\sum_{i=N+1}^{N+M}P_{i-N}(\bm{x}_j) \dot{u}_i(t)=-\left(\mathcal{L}\hat u\right)(t,\bm{x}_j) + f(t,\bm{x}_j) \label{coll}
\end{equation}
for $j\in \mcI_{\partial}$
\begin{equation}
     \displaystyle 0=-\left(\mathcal{B}\hat u\right)(t,\bm{x}_j) + g(t,\bm{x}_j), \label{bcond}
\end{equation}
Moreover, by imposing the orthogonality conditions, we have,
for $j=N+1,N+2,\dots, N+M$
\begin{equation}
      0= \sum_{i=1}^NP_{j-N}({\bm x}_i) u_i(t).\label{ort}
\end{equation}

We note that if the boundary conditions change on different parts of $\partial^*\Omega$ then we have one equation of type \eqref{bcond} for each part; equations \eqref{bcond} and \eqref{ort} depend only on $\bm u$.

In compact notation, we can write 

\begin{equation}
    \label{eq_matrice}
    A\bm{\dot{\bm u}}(t)=\bm{b}(t,\bm u),
\end{equation}
%$$ \hat u(t,\bm{x})=\displaystyle\sum_{i=1}^Nu_i(t)\phi(\norm{\bm{x}-\bm{x_i}})+\displaystyle\sum_{i=N+1}^{N+M}u_i(t)P_{i-N}(\bm{x}),\qquad (t,\bm x)\in \rr_0^+\times \overline \Omega$$
where $A\in\rr^{(N+M)\times(N+M)}$ has entries

\begin{equation}
a_{j,i}=\left\{
\begin{array}{ll}
    \phi\left(\norm{\bm x_j-\bm x_i}\right), & 1\le i\le N,\ j\in \mcI_0,\\  P_{i-N}(\bm x_j) , & N+1\le i\le N+M,\ j\in \mcI_0,\\
    0, & \mbox{otherwise},
\end{array}
\right.\label{CoeA}
\end{equation}

and  $\bm b(t,\bm u)\in\rr^{(N+M)}$ has entries

\begin{equation}
b_{j}(t,\bm u)=\left\{
\begin{array}{ll}
    -(\mcL \hat u)(t,\bm x_j) + f(t,\bm x_j), & j\in \mcI_0,\\   -(\mathcal{B} \hat u)(t,\bm x_j) + g(t,\bm x_j), & j\in \mcI_\partial,\\
    \sum_{i=1}^NP_{j-N}({\bm x}_i) u_i(t), & N+1\le j\le N+M.
\end{array}
\right.\label{Coeb}
\end{equation}

\section{Problem Approximation}
\label{sec_PA}
\begin{figure}
\centering
\begin{tabular}{cc}
\includegraphics[width=.45\textwidth]{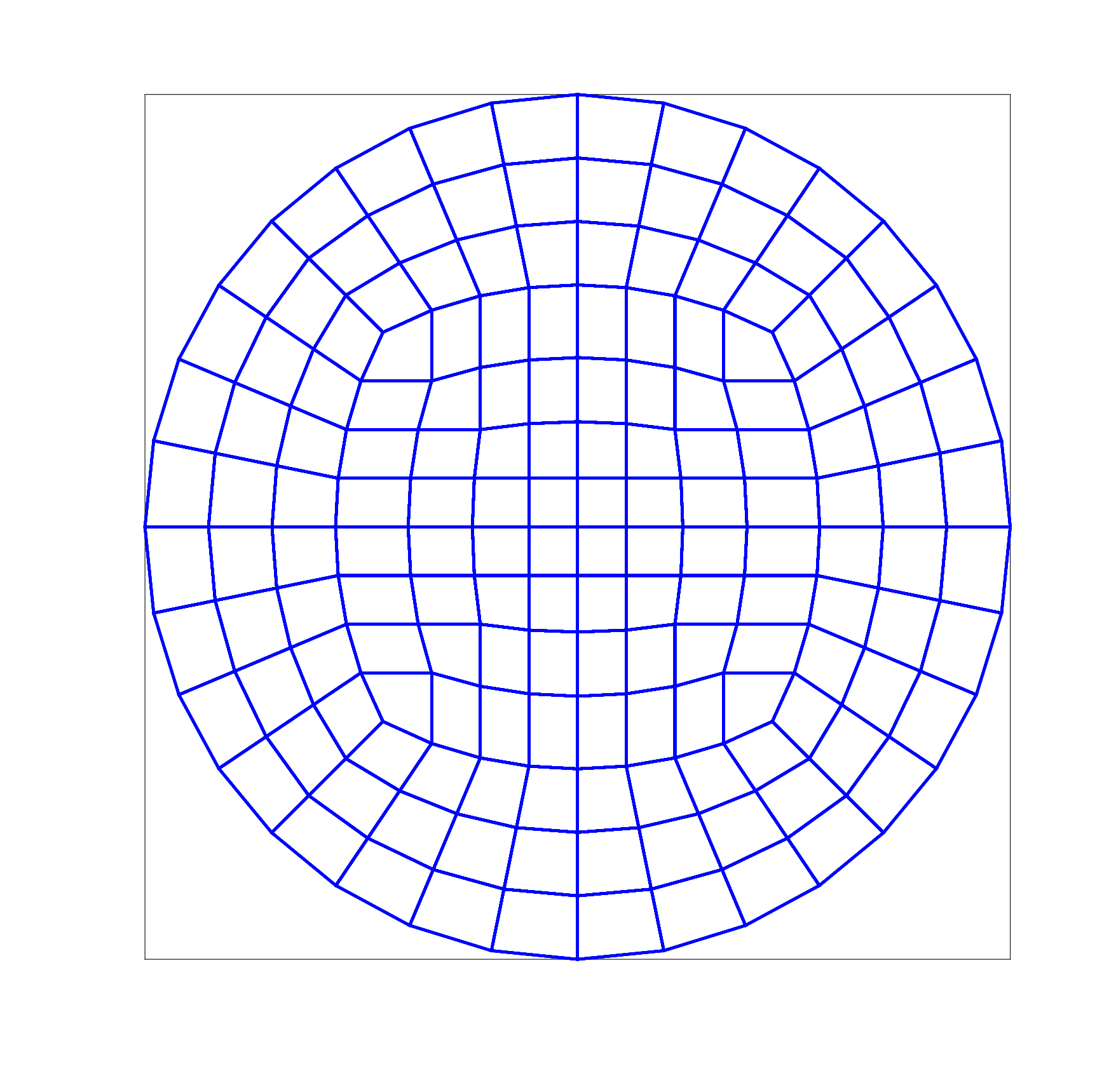} &
\includegraphics[width=.45\textwidth]{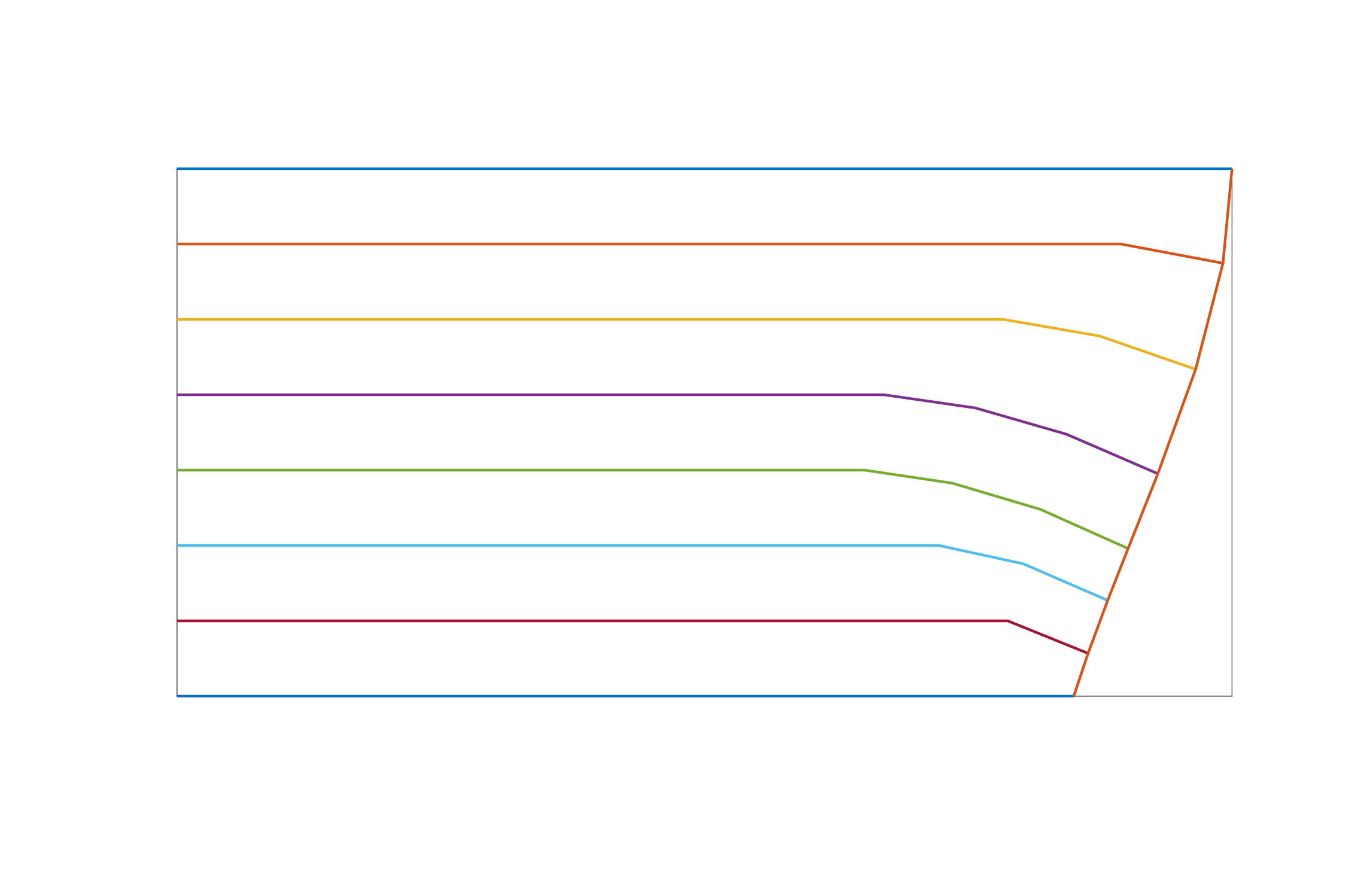}\\
    (a) & (b)
\end{tabular}
\caption{Upper face discretization (a). Cross section view and curves where the horizontal discretization is repeated (b).}
\label{fig_griglia}
\end{figure}

Generally, the filter is a geometrical object with cylindrical symmetry, having circular horizontal sections and a lateral profile dictated by the type of filter; this profile can be linear or curvilinear. Hence, we can treat the coffee pod as a spatial domain $\mathcal{D}$ having cylindrical symmetry, like the one given in Figure~\ref{CofPod}, and its shape may vary depending on the brand. In order to apply the collocation methods, we need to define the collocation nodes $\mcX$. The nodes may be generated by using an analytical approach that creates and optimizes an exaedral mesh in a spatial domain that reproduces a filter.  On the upper face of the filter, which is circular, the collocation nodes are taken as follows: in the center there is a square pattern, after a certain distance from the center of the circumference,  the nodes are taken on arcs of ellipses with increasing curvature as they move away from the center, until they approach the circular boundary of the domain as shown in Figure \ref{fig_griglia}(a). On the lower face, which also has a circular shape, the collocation nodes are taken so that they have the same pattern as those on the upper face; the nodes on the two faces are in one-to-one correspondence.

Starting from the upper face, all points are projected onto the longitudinal section below, and then coordinates $x_1$ and $x_2$ are appropriately rescaled. These points are not coplanar, i.e., they do not belong to the same plane orthogonal to $x_3$, since the heights of points are calculated so that the direction of the segment linking the point on the boundary and the previous inner point is orthogonal to the boundary, as shown in Figure \ref{fig_griglia}(b).

\begin{figure}{H}
\centering
\begin{tabular}{cc}
    \includegraphics[width=0.49\textwidth]{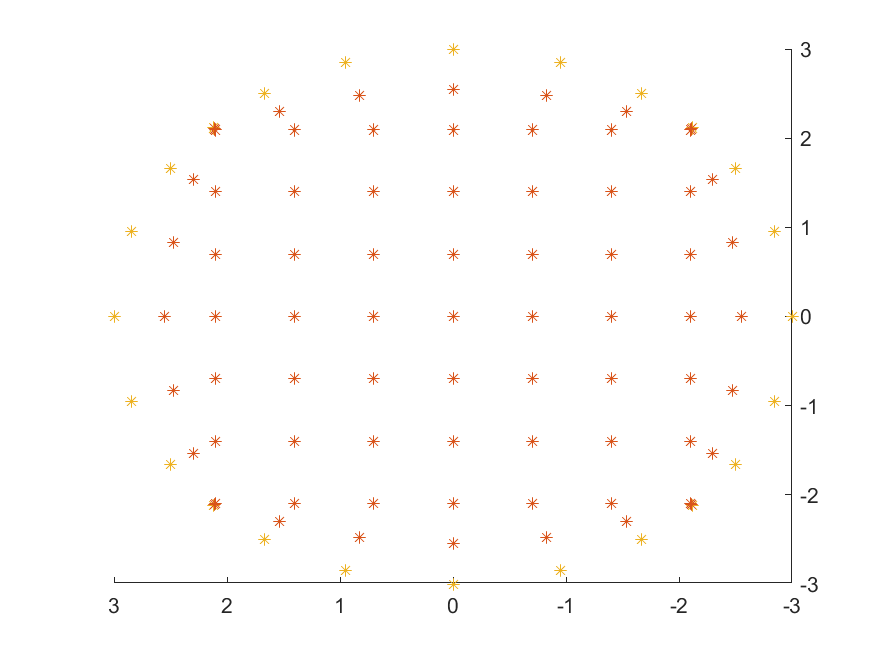} &
	\includegraphics[width=0.49\textwidth]{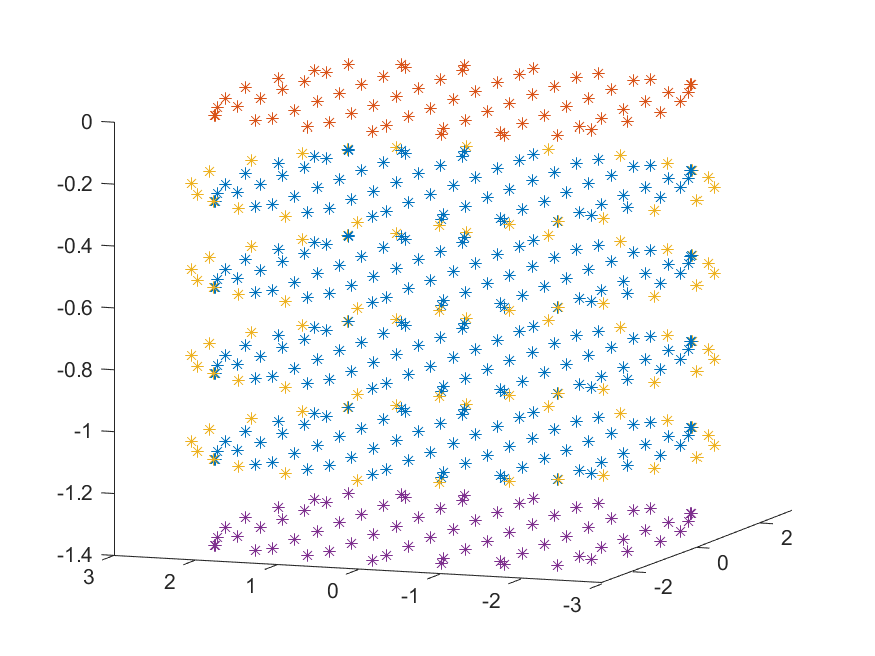}\\
    (a) & (b)
\end{tabular}
	\caption{Collocation nodes of the Percolation Model when $\mathcal{D}$ is a cylinder with height $L=1.388$ and $R=3$.}
    \label{figure_nodes_filtro}
\end{figure}

In this first work, we consider a cylindrical filter as in Figure \ref{CofPod}, and we consider a Cartesian coordinate system centered at the center of the upper face $\Gamma_1$ of $\mathcal{D}$ and with $x_3$-axis pointing upwards. Let $\Gamma_2$ and $\Gamma_3$ be the lateral and the lower surface of $\mathcal{D}$ and $\mathcal{C}$ the interior of $\mathcal{D}$.   
To fix the ideas, if $\mathcal{D}$ is a circular cylinder of height $L$ and radius $R$, Figure \ref{figure_nodes_filtro} shows the discretization of the domain. 
In detail,
Figure \ref{figure_nodes_filtro} shows nodes obtained with 6 slices, each with 97 collocation nodes in a cylindrical domain with $L=1.388cm$ and $R=3cm$.

Then, with the notation introduced in the previous section, we have that $\Omega=\mathcal{C}$, $\partial\Omega=\Gamma_1\cup \Gamma_2\cup \Gamma_3$, $\overline\Omega=\mathcal{D},$ $\mcX$ is the set of the vertices of the exaedrals in the mesh above described, $N$ is the cardinality of $\mcX$, and  $\mcI_\partial=\mcI_1\cup\mcI_2\cup\mcI_3 $ with $\mcI_k=\{i\in\mcI\ :\  \bm x_i\in\mcX\cap \Gamma_k\},$ $k=1,2,3.$ 
Hence, $\mcI_0$ contains the indices of the nodes of $\mcX$ in $\mathcal{C}$, and for $k=1,2,3,$ $\mcI_k$ contains the indices of the nodes of $\mcX$ in the boundary $\Gamma_k$. Moreover we have that $\{\mcI_0,\mcI_1,\mcI_2,\mcI_3\}$ is a partition of $\mcI=\{1,2,\dots,N\},$ and we define $\mcI_P=\{N+1,N+2,\dots, N+M\}.$ 
In the example shown in Figure \ref{figure_nodes_filtro}, we have $\vert{\mathcal{I}_0}\vert=292$, $\vert{\mathcal{I}_1}\vert=73, \vert{\mathcal{I}_2}\vert=96, \vert{\mathcal{I}_3}\vert=73$ the total number of collocation nodes is $\vert{\mathcal{I}_0}\vert+ \vert{\mathcal{I}_1}\vert+ \vert{\mathcal{I}_2}\vert+ \vert{\mathcal{I}_3}\vert=534$.
To avoid conflicting boundary conditions, we have no nodes in  $\Gamma_1\cap\Gamma_2$ and $\Gamma_2\cap\Gamma_3.$ 

First of all, we give the approximation method obtained by applying Kansa's method to our considered problem. 

By considering the first two equations in \eqref{model} with boundary conditions in \eqref{RicDar}, we have a problem of kind \eqref{PDE}, then we denote the corresponding differential operators with $\mcL_h$, the three boundary differential operators with $\mcB_{h1},$ $\mcB_{h2},$ and $\mcB_{h3}$, with corresponding boundary indices in $\mcI_{1},$ $\mcI_{2},$ and $\mcI_{3}$, and given functions $f_h$, $g_{h1},$ $g_{h2},$ and $g_{h3}.$ By approximating $h$ as in \eqref{hatu}, let be $\bm h:\rr_0^+\to \rr^{N+M}$ then the corresponding system of DAEs is

\begin{equation}
    \label{eq_matrice_problema_test}
    A\bm{\dot{\bm h}}(t)=\bm{b}_h(t,\bm h),
\end{equation}
where $A\in\rr^{(N+M)\times(N+M)}$ is given in \eqref{CoeA}
and  $\bm b_h(t,\bm h)\in\rr^{(N+M)}$ has entries

\begin{equation}
b_{j}(t,\bm u)=\left\{
\begin{array}{ll}
    -(\mcL_h \hat h)(t,\bm x_j) + f_h(t,\bm x_j), & j\in \mcI_0,\\   -(\mathcal{B}_{h1} \hat h)(t,\bm x_j) + g_{h1}(t,\bm x_j), & j\in \mcI_1,\\
    -(\mathcal{B}_{h2} \hat h)(t,\bm x_j) + g_{h2}(t,\bm x_j), & j\in \mcI_2,\\
    -(\mathcal{B}_{h3} \hat h)(t,\bm x_j) + g_{h3}(t,\bm x_j), & j\in \mcI_3,\\
    \sum_{i=1}^NP_{j-N}({\bm x}_i) u_i(t), & N+1\le j\le N+M.
\end{array}
\right.\label{Coebh}
\end{equation}

In the same way, we discretize the other equations in \eqref{model}, obtaining the complete system of DAEs.

The principal parameters of the model are the pressure and temperature of incoming water, the granulometry of the coffee powder, and the type of coffee blend; some other parameters depend on those. In the algorithm implementation, we have considered Arabica coffee and “optimal” granulometry, and we studied the behavior of the caffeine in the solid and liquid phases. In the Table \ref{tab_valori}, the values of the variables are reported; moreover, since we are referring to a particular chemical species, we had to replace the superscripts and the subscripts $k$ with $1$.

\begin{table}[]
    \centering
    \begin{tabular}{lcr|lcr}
    \toprule
    variable & value & unit & variable & value & unit\\
    \midrule
        $p_{z0}$ & $6$ & $b$ &
        $T_{z0}$ & $88$ & ${^\circ} C$\\
        $\epsilon$ & $0.305$ & &
        $k$ & $1.8282$ & $cm/d$ \\
        $\overline{t}$ & $2.3148\cdot 10^{-04}$ & $d$&
        $\rho_0$ & $0.01$ & $Kg/m^3$\\
        $h_{z0}$ & $6118.3$ & $cm$&
        $a_1$ & $3184.9$ & $1/d$\\
        $S_0$ & $10^{-5}$ & $1/cm$&
        $\beta_T^1$ & $10$ & $cm$\\
        $\beta_L^1$ & $100$ & $cm$&
        $D_1$ & $86400\cdot10^{-5}$ & $cm^2/d$\\
        $\rho c$ & $4.18\cdot10^{-3}$ & $J/cm^3 {^\circ} C$&
        $\rho^s c^s$ & $3.184\cdot 10^{-3}$ & $J/cm^3{^\circ} C$\\
        $\lambda$ & $86400\cdot5\cdot10^{-3}$ & $J/d\,cm{^\circ} C$&
        $f_\mu$ & $1$ & $ $\\
        $\chi$ & $0$ & $ $&
        $\Phi_h$ & $86400\cdot6.5\cdot10^{-5}$ & $1/d$\\
        $\Phi_1$ & $259200 $ & $cm/d$&
        $T_0$ & $70$ & $ {^\circ} C$\\
        $C_{10}^s$ & $0.01254$ & $Kg/L$& \\
    \bottomrule
    \end{tabular}
    \caption{Variable values of the problem for caffeine.}
    \label{tab_valori}
\end{table}

\section{Numerical results}
\label{sec_NR}
The results obtained were compared with those obtained by solving the model using the FEFLOW software~\cite{feflow_sito}, which we consider correct since they are supported by the values obtained from laboratory chemical analyses~\cite{Art_perc}. In the algorithm, we made two changes with respect to the model: the initial hydraulic head profile was considered constant and equal to $h_{z0}$, and the initial condition for the temperature in the $\Gamma_1$ face was changed into $T_{z0}$ to ensure that the compatibility conditions were satisfied.

\FloatBarrier
\bigskip
\textbf{Hydraulic Head}

\begin{figure}
    \centering
    \subfloat
    {\includegraphics[width=0.495\linewidth]{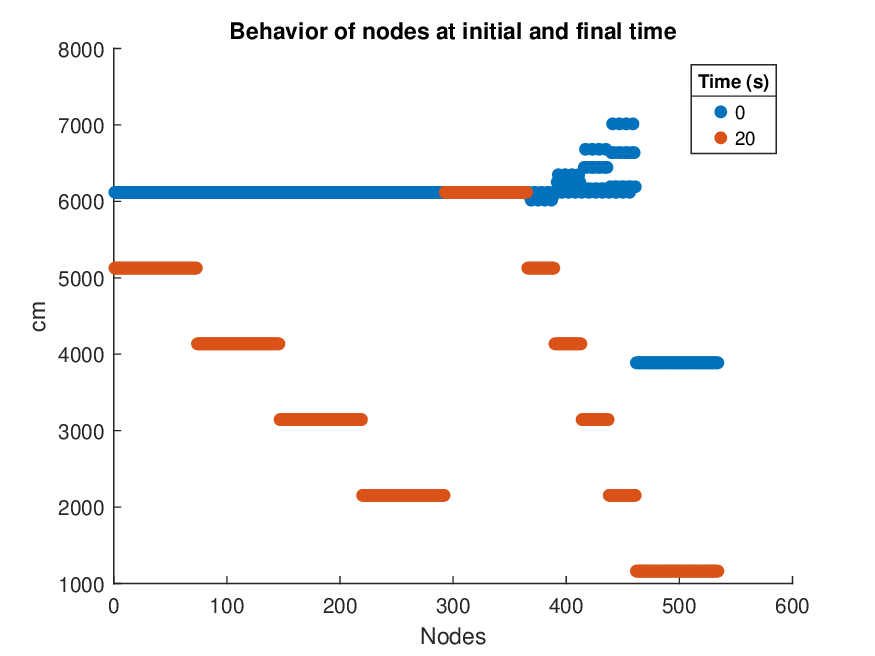}}
    \subfloat
    {\includegraphics[width=0.495\linewidth]{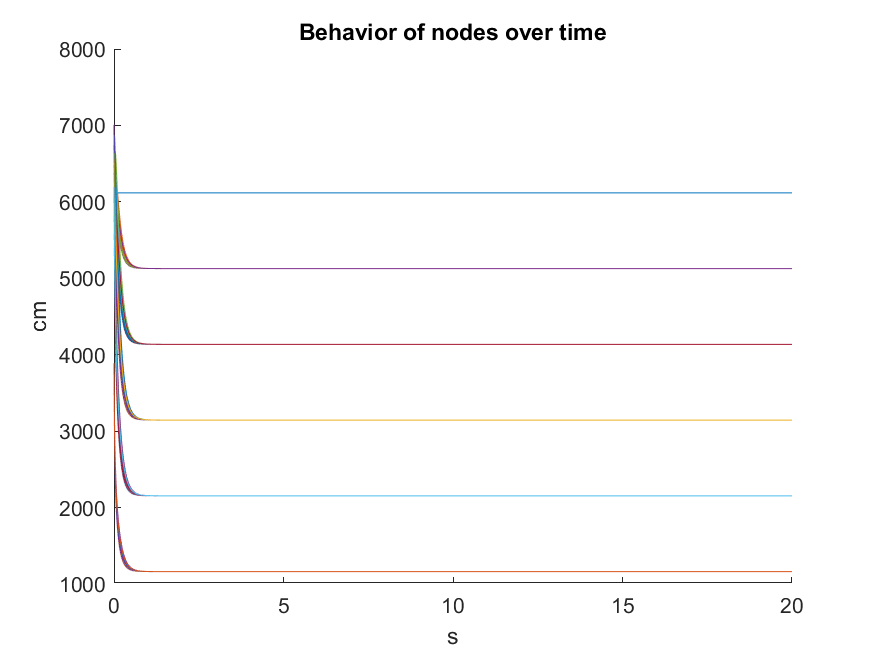}}
    \caption{Trend of the approximate Hydraulic Head over time.}
    \label{fig_h_su_tempo}
\end{figure}

In Figure \ref{fig_h_su_tempo}, we can observe the evolution of the approximate hydraulic head over time. 

The left side represents the value of the approximated hydraulic head in each node at the initial (blue) and final (orange) time of the percolation process. As can be noticed, at the beginning of the process, the value of the hydraulic head at some nodes is not exactly equal to the initial condition we imposed in the model. This occurs because MATLAB automatically adjusts the initial conditions in order to solve the system. On the right side of Figure~\ref{fig_h_su_tempo}, we can see the evolution of hydraulic head in each collocation node over time. From these two figures, it is possible to see that the nodes at the same height behave in the same way. It can also be observed that, within less than one second, the hydraulic head inside the filter decreases and then remains constant until the end of the process. In Table \ref{tab_errore_h}, a quantitative comparison of the results at the end of the percolation is shown. To do this, since nodes at the same height show exactly the same value for the hydraulic head, we can compare the results obtained in MATLAB with those in FEFLOW by considering the values in two nodes at the same height.

\begin{table}[]
    \centering
    \begin{tabular}{cccc}
    \toprule
    Node Height & FEFLOW & MATLAB & Error\\
    (cm) & (cm) & (cm) & (\%)\\
    \midrule
    0 & 6118.39 & 6118.29 & 0.00016\\
    -0.2776 & 5125.76 & 5127.11 & 0.0263\\
    -0.5552 & 4133.13 & 4135.92 & 0.0675\\
    -0.8328 & 3140.5 & 3144.73 & 0.13\\
    -1.1104 & 2147.87 & 2153.54 & 0.26\\
    -1.388 & 1155.24 & 1162.36 & 0.62\\
    \bottomrule
    \end{tabular}
    \caption{Relative Error of Hydraulic Head.}
    \label{tab_errore_h}
\end{table}

\FloatBarrier
\textbf{Solid Caffeine}

\begin{figure}
    \centering
    \subfloat
    {\includegraphics[width=0.495\linewidth]{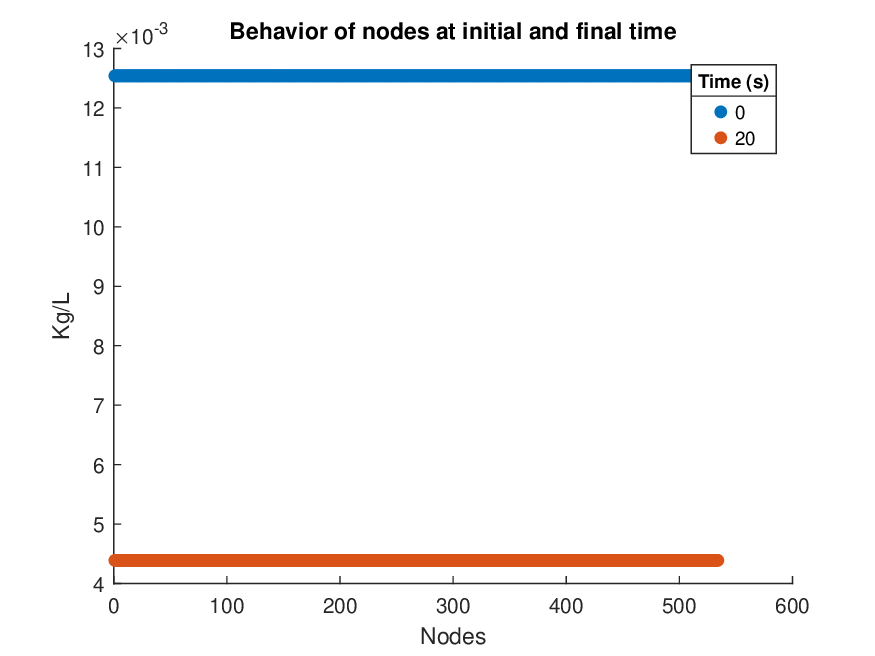}}
    \subfloat
    {\includegraphics[width=0.495\linewidth]{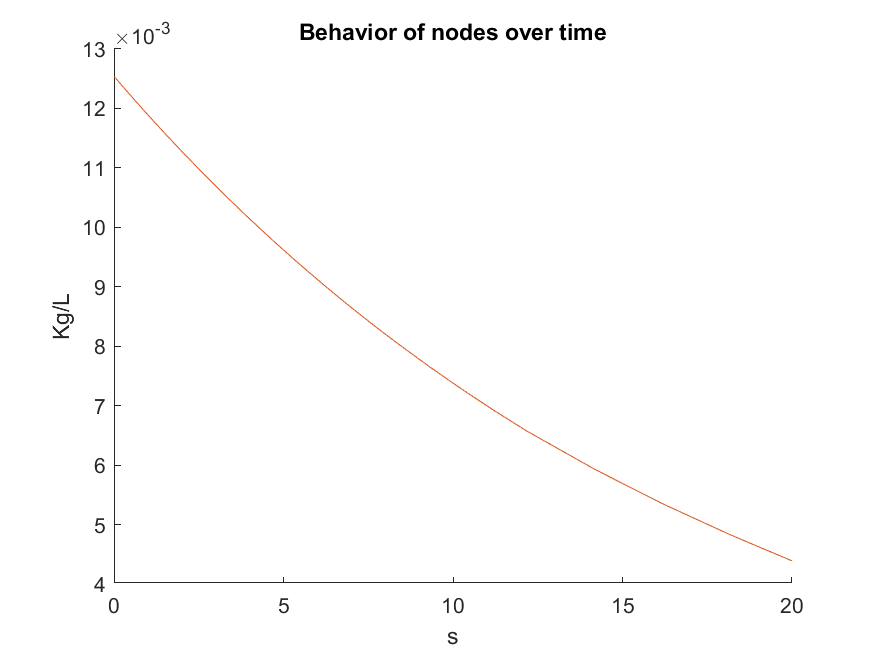}}
    \caption{Trend of the approximate Solid Caffeine over time.}
    \label{fig_cf_su_tempo}
\end{figure}

In Figure \ref{fig_cf_su_tempo}, we can see the trend of the approximate solid caffeine over time. As before, the left side shows the value of solid caffeine in each node at the beginning (blue) and end (orange) of the percolation process, and on the right side, we can observe the evolution of solid caffeine in each collocation node over time. As can be observed, all nodes have exactly the same behavior throughout the entire domain. The value obtained from MATLAB at the end of the process is $4.3888 \cdot 10^{-3} Kg/L$, while the value obtained from FEFLOW is $4.32814\cdot 10^{-3} Kg/L$, thus we have an error of $1.4\%$.

\begin{figure}
    \centering
    \subfloat
    {\includegraphics[width=0.495\linewidth]{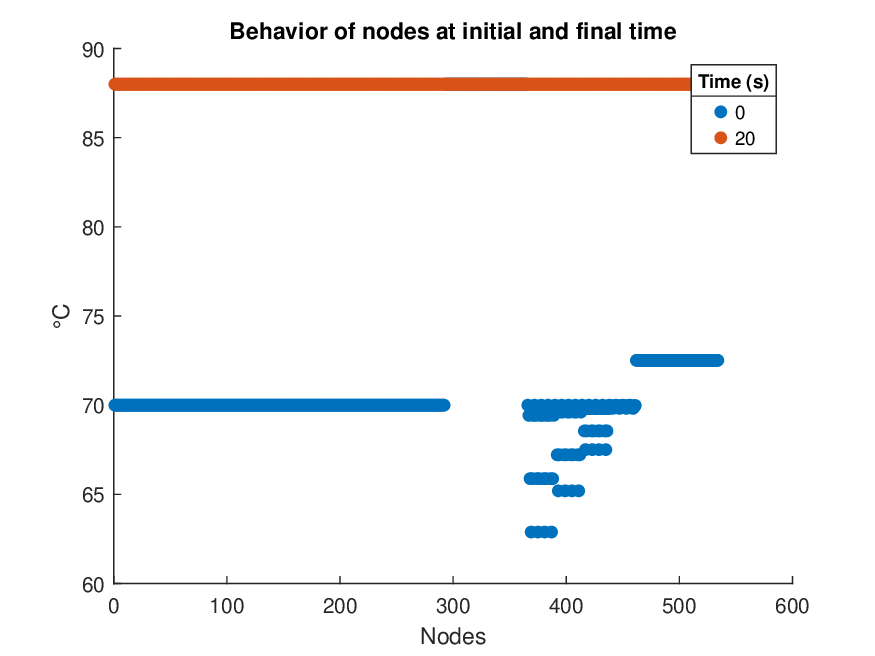}}
    \subfloat
    {\includegraphics[width=0.495\linewidth]{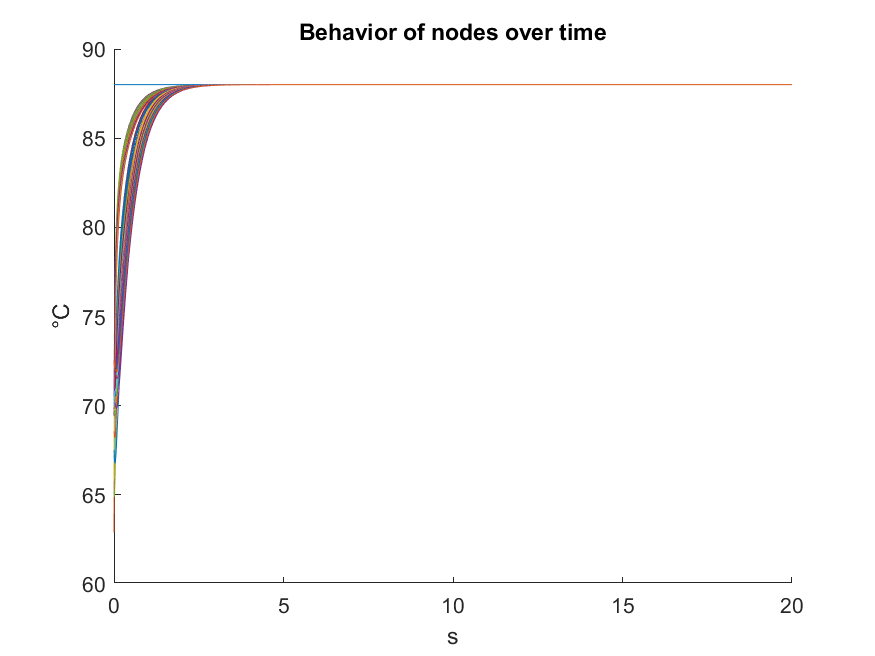}}
    \caption{Trend of the approximate Temperature over time.}
    \label{fig_T_su_tempo}
\end{figure}

\FloatBarrier
\textbf{Temperature}\\
In Figure \ref{fig_T_su_tempo}, we can observe the evolution of the approximate temperature.

The figure on the left side shows the values of temperature in each node at the beginning (blue) and end (orange) of the process. As in the hydraulic head solution, some initial values differ from the original ones because MATLAB has modified them. In the right side figure, the evolution of temperature in each collocation node over time is shown, and it can be seen how, in the first three seconds, the temperature of the coffee pods equalizes with that of the incoming water, after which it remains constant until the end of the process. A comparison of the results obtained in MATLAB and FEFLOW at the end of the percolation is shown in Table \ref{tab_errore_T}.

\begin{table}[]
    \centering
    \begin{tabular}{cccc}
    \toprule
    Node Height & FEFLOW & MATLAB & Error\\
    (cm) & ({$^\circ C$}) & ({$^\circ C$}) & (\%)\\
    \midrule
    0 & 88 & 88 & 0\\
    -0.2776 & 87.933 & 88 & 0.0761\\
    -0.5552 & 87.866 & 88 & 0.15\\
    -0.8328 & 87.798 & 88 & 0.23\\
    -1.1104 & 87.73 & 88 & 0.31\\
    -1.388 & 87.66 & 88 & 0.39 \\
    \bottomrule
    \end{tabular}
    \caption{Relative Error of Temperature.}
    \label{tab_errore_T}
\end{table}

\FloatBarrier
\textbf{Liquid Caffeine}

\begin{figure}
    \centering
    \subfloat
    {\includegraphics[width=0.495\linewidth]{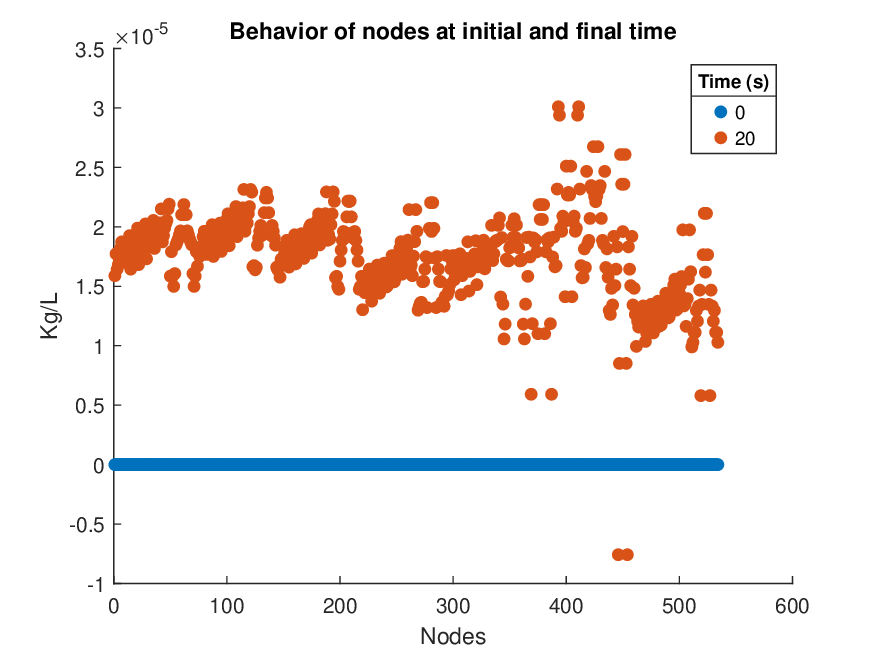}}
    \subfloat
    {\includegraphics[width=0.495\linewidth]{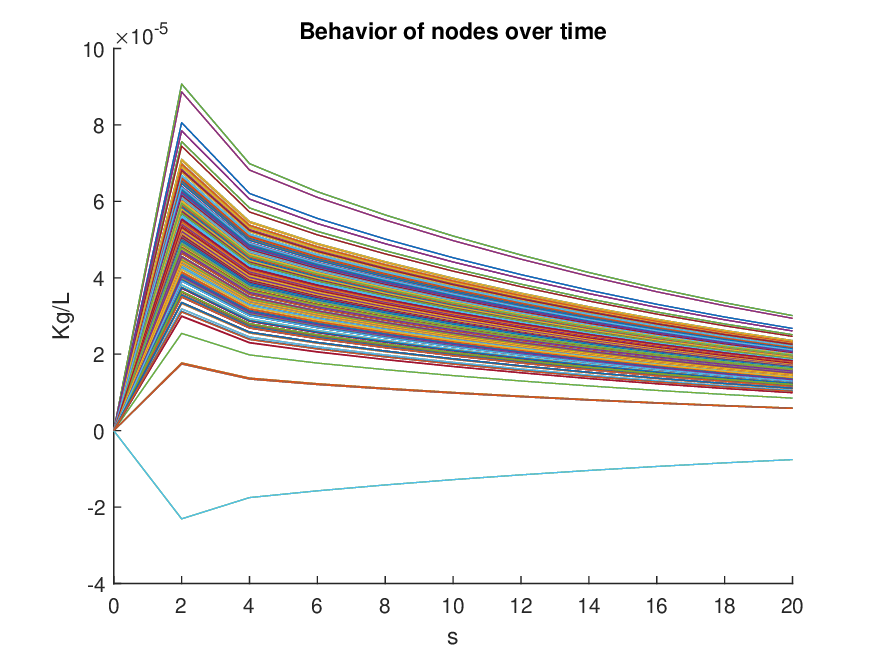}}
    \caption{Trend of the approximate Liquid Caffeine over time.}
    \label{fig_cf_liquida_su_tempo}
\end{figure}

The evolution of the liquid caffeine at the beginning, during, and at the end of the percolation process is shown in Figure \ref{fig_cf_liquida_su_tempo}. In particular, on the left side, the blue profile gives the value of the concentration at each node at the initial time, while the orange profile gives the value of the concentration at each node at the final time. On the right side, the behavior of the liquid caffeine in each node over time is reported.

\begin{figure}
    \centering
    \includegraphics[width=0.5\linewidth]{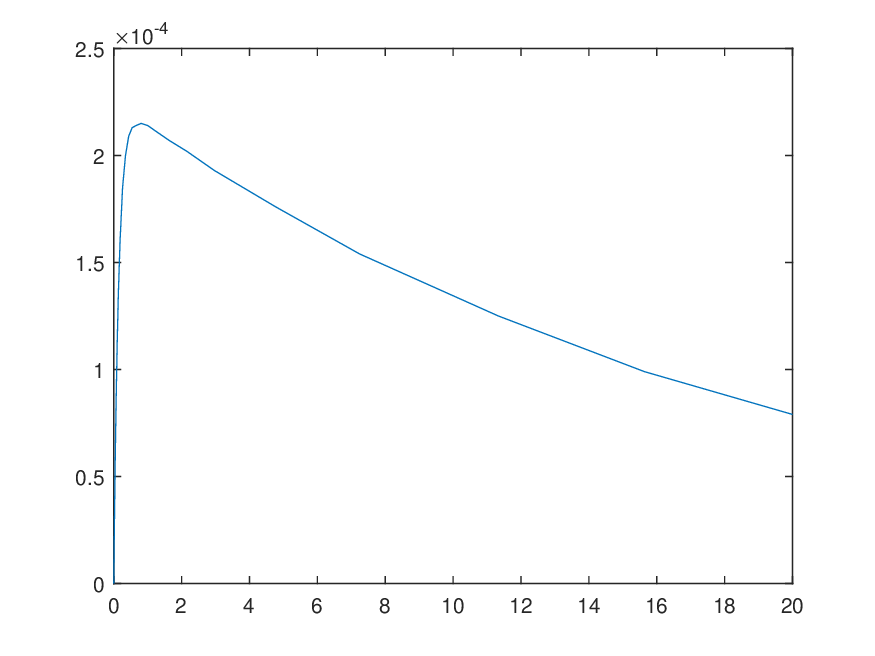}
    \caption{Trend of Liquid Caffeine.}
    \label{fig_cf_liquida_vera}
\end{figure}

Looking at the left side of the figure, negative values appear in two nodes. Naturally, this is not physically possible, and the phenomenon is most likely due to numerical errors. From the results of FEFLOW, we know that the expected trend and value of the liquid caffeine in a node over time is as illustrated in Figure \ref{fig_cf_liquida_vera}, with small variations in the value depending on the height of the node.

If we compare the behavior of approximate liquid caffeine in each collocation node over time, i.e., Figure \ref{fig_cf_liquida_su_tempo}, and the expected value, in Figure \ref{fig_cf_liquida_vera}, we can see how the curves show trends that are in good agreement, except for the initial slope and the magnitude of the liquid caffeine, which in our results is the half or, in some nodes, up to one order less.

\section{Conclusion}
\label{sec_conclu}
In this paper, we considered the problem of coffee percolation. We developed a model
that describes the process from a physico-chemical perspective. Specifically, we analyzed
the behavior of the hydraulic head, the evolution of chemicals in solid and liquid
phases, and temperature dynamics. The collocation method with RBFs,
is used to obtain a discretization scheme for the percolation problem. We implemented
this algorithm in MATLAB to compute the numerical solution of the problem. The results were
compared with those obtained by using FEFLOW, whose results are considered reliable
as they had been validated through laboratory chemical analysis. The results obtained for hydraulic head, solid caffeine, and temperature are extremely good, while those for liquid caffeine are somewhat less accurate and must be refined in future work.
A possible future development is the refinement of the approximation for liquid substances, since, given the very good agreement for the other results, we expect to be
able to improve this as well. Afterward, other chemical substances of interest could be
incorporated into the system. Another interesting study is the implementation of other
RBF families to compare the results in terms of stability and accuracy. Moreover, it is
possible to solve the model using other methods, such as the Finite Element Method, to
compare the results. Furthermore, these results could be compared with those obtained
from other models, such as the one proposed by Moroney~\cite{Moroney}, or Cameron~\cite{cameron2020systematically}. 
The model can be further improved by considering the interaction between chemical
species and the transport of fine particles that change the porosity of the medium, thus
creating the compact layer in the coffee pod. It is worth mentioning that the proposed
percolation model can be easily generalized to a wide range of flow problems through
porous media, finding application in many areas, such as hydrogeological fields.
Finally, it can be observed that this model represents an ad hoc tool that opens the way
for the customization of espresso coffee taste, which is very important in the industrial
context.
\bibliographystyle{unsrt} %il comando unsrt ordina le citazioni in ordine di citazione
\bibliography{Tesi_magistrale_Pacini_Gianluca/BIBLIO}

\end{document}